\RequirePackage{ifpdf}
\ifpdf 
\documentclass[pdftex]{sigma}
\else
\documentclass{sigma}
\fi

\begin{document}

\allowdisplaybreaks

\renewcommand{\thefootnote}{$\star$}

\renewcommand{\PaperNumber}{067}

\FirstPageHeading

\ShortArticleName{Symplectic Applicability of Lagrangian Surfaces}

\ArticleName{Symplectic Applicability of Lagrangian Surfaces\footnote{This paper is a
contribution to the Special Issue ``\'Elie Cartan and Dif\/ferential Geometry''. The
full collection is available at
\href{http://www.emis.de/journals/SIGMA/Cartan.html}{http://www.emis.de/journals/SIGMA/Cartan.html}}}

\Author{Emilio~MUSSO~$^\dag$ and Lorenzo NICOLODI~$^\ddag$}

\AuthorNameForHeading{E Musso and L. Nicolodi}

\Address{$^\dag$~Dipartimento di Matematica, Politecnico di Torino,\\
\hphantom{$^\dag$}~Corso Duca degli Abruzzi 24, I-10129 Torino, Italy}
\EmailD{\href{mailto:emilio.musso@polito.it}{emilio.musso@polito.it}}

\Address{$^\ddag$~Di\-par\-ti\-men\-to di Ma\-te\-ma\-ti\-ca,
Uni\-ver\-si\-t\`a degli Studi di Parma,\\
\hphantom{$^\ddag$}~Viale G.P.\ Usberti 53/A,
I-43100 Parma, Italy}
\EmailD{\href{mailto:lorenzo.nicolodi@unipr.it}{lorenzo.nicolodi@unipr.it}}

\ArticleDates{Received February 25, 2009, in f\/inal form June 15, 2009;  Published online June 30, 2009}

\Abstract{We develop an approach to af\/f\/ine symplectic invariant
geometry of Lagrangian surfaces by the method of moving frames. The
fundamental invariants of elliptic Lagrangian immersions in af\/f\/ine
symplectic four-space are derived together with their integrability
equations. The invariant setup is applied to discuss the question of
symplectic applicability for elliptic Lagrangian immersions.
Explicit examples are considered.}

\Keywords{Lagrangian surfaces; af\/f\/ine symplectic geometry; moving
frames; dif\/ferential invariants; applicability}

\Classification{53A07; 53B99; 53D12; 53A15}

\renewcommand{\thefootnote}{\arabic{footnote}}
\setcounter{footnote}{0}

\section{Introduction}

\looseness=1
While the subject of metric invariants of Lagrangian submanifolds
has received much attention in the literature, the subject of af\/f\/ine
symplectic invariants of Lagrangian submanifolds has been less
studied (cf.~\cite{MK} for recent work in this direction). Early
contributions to symplectic invariant geometry of submanifolds go
back to Chern and Wang \cite{Chsymplectic}, who studied submanifolds
in projective $(2n+1)$-space under the linear symplectic group. Some
recent works related to the study of af\/f\/ine symplectic invariants,
mainly for the case of curves and hypersurfaces in Euclidean space,
include \cite{De, KOT, AD}. The aim of
this paper is to develop an approach to af\/f\/ine symplectic invariant
geometry of Lagrangian surfaces in standard af\/f\/ine symplectic
four-space, using the method of moving frames. The invariant setup
is then applied to discuss the problem of symplectic applicability
for Lagrangian immersions.

\looseness=1
Let $M$ be a 2-dimensional manifold and let $f : M\to \mathbb{R}^4$ be a
Lagrangian immersion into standard af\/f\/ine symplectic four-space
$(\mathbb{R}^4,\Omega)$. The basic observation is that the natural Gauss map
of $f$, whose value at a point is determined by the 1-jet of $f$ at
that point, takes values in the oriented Lagrangian Grassmannian of
$(\mathbb{R}^4, \Omega)$, the manifold $\Lambda_2^+$ of oriented Lagrangian
2-subspaces in $\mathbb{R}^4$. The oriented Lagrangian Grassmannian admits
an alternative description as the conformal compactif\/ication of
Minkowski 3-space $\mathbb{R}^{2,1}$, the projectivization of the positive
nullcone of $\mathbb{R}^{3,2}$. The local isomorphisms between the group
$\mathrm{Sp}(4,\mathbb{R})$ of linear symplectomorphisms and the conformal Lorentz
group $\mathrm{O}(3,2)$ establishes a close relationship between the
symplectic geometry of a Lagrangian surface in $\mathbb{R}^4$ and the
conformal Lorentzian geometry of its Gauss map. Accordingly, a
Lagrangian immersion is of \textit{general type} (resp.,
\textit{special type}) if its  Gauss map is an immersion (resp., has
constant rank one). Moreover, a Lagrangian immersion of general type
is {\it elliptic} (resp., {\it hyperbolic}, {\it parabolic}) if its
Gauss map is a spacelike (resp., timelike, lightlike) immersion into
$\Lambda_2^+$.

In this work we study elliptic Lagrangian immersions
in $(\mathbb{R}^4,\Omega)$ under the group $\mathbb{R}^4 \rtimes \mathrm{Sp}(4,\mathbb{R})$
of af\/f\/ine symplectic transformations.
In Section~\ref{s:pre}, we collect the background material on af\/f\/ine
symplectic frames and their structure equations, and brief\/ly
describe the conformal geometry of the Lagrangian Grassmannian~$\Lambda_2^+$.

\looseness=1
In Section \ref{s:reduction}, we develop the moving frame method for
elliptic Lagrangian immersions in af\/f\/ine symplectic geometry. If~$f
: M \to \mathbb{R}^4$ is an elliptic Lagrangian immersion, we consider on
$M$ the complex structure def\/ined by the conformal structure induced
by the Gauss map of~$f$ and a choice of orientation on~$M$. After
the third frame reduction, we make a successive reduction with
respect to a f\/ixed local complex coordinate in the Riemann surface
$M$, which yields a unique adapted frame f\/ield along~$f$. From this
we derive the three local af\/f\/ine symplectic invariants for $f$,
namely the complex-valued smooth functions $t$, $h$, $p : M \to \mathbb{C}$,
and establish their integrability equations, which give the
existence and uniqueness theorem for elliptic Lagrangian immersion
in af\/f\/ine symplectic geometry. At this stage, a natural question to
ask is to what extent the invariants are actually needed to
determine the elliptic Lagrangian surface up to symmetry (cf.\
Bonnet's problem in Euclidean geometry). We will look at this
problem from the classical point of view of applicable surfaces,
a concept generalized by \'E. Cartan to $G$-deformations of submanifolds
of any homogeneous space $G/G_0$ (cf.~\cite{Fu,Ca1,Ca2,Gr, JensenJDG}).

\looseness=1
In Section \ref{s:applicability}, we apply the invariant setup to
discuss rigidity and applicability of elliptic Lagrangian
immersions. Already after the second reduction, we associate with
any elliptic Lag\-rangian immersion $f$ a cubic dif\/ferential form
$\mathcal{F}$, the \textit{Fubini cubic form}. Proposition
\ref{prop:gen} shows that the position of a generic elliptic
Lagrangian immersion is completely determined by $\mathcal{F}$, up
to af\/f\/ine symplectic transformation. The elliptic Lagrangian
immersions which are not determined by the Fubini cubic form alone
are called \textit{applicable}. Theorem \ref{thm:applicable} relates
applicability to the complex structure of $M$, using the normalized
Hopf dif\/ferential $\mathcal{H}$, a quadratic dif\/ferential form
naturally associated to any elliptic Lagrangian immersion with
respect to third order frames. More precisely, it proves that an
elliptic Lagrangian immersion is applicable if and only if there
exists on~$M$ a non-zero holomorphic quadratic dif\/ferential
$\mathcal{Q}$ such that~$\mathcal{H}$ and~$\mathcal{Q}$ are linearly
dependent over the reals at each point of~$M$. An interesting
consequence of this result is that applicable elliptic Lagrangian
immersions admit an associated family of noncongruent elliptic
Lagrangian immersions parametrized by~$\mathbb{R}$, which in turn suggests
their integrable character. This is not surprising when one
considers that the Gauss map of an applicable elliptic Lagrangian
immersion is an isothermic spacelike immersion of~$\Lambda_2^+$.
Indeed, it is well known that a spacelike isothermic immersion in
3-dimensional conformal space comes in an asso\-cia\-ted fa\-mi\-ly of
(second order) deformations which coincide with the classical
$T$-transforms of Bianchi and Calapso (cf.~\cite{Mtrieste,HMN,MNhouston}). The equations governing elliptic
Lagrangian immersions are also derived. We then consider the class
of elliptic Lagrangian surfaces with totally umbilic Gauss map. They
are def\/ined by the vanishing of the normalized Hopf dif\/ferential~$\mathcal{H}$. Proposition~\ref{prop:non-congruentTU} shows that, up
to a biholomorphism of $M$, two noncongruent, totally umbilic,
elliptic Lagrangian immersions are applicable. Finally, Proposition~\ref{prop:TU-complexcurve} proves that a complex curve in~$\mathbb{C}^2$
without f\/lex points is a totally umbilic elliptic Lagrangian
immersion and that, conversely, any totally umbilic elliptic
Lagrangian immersion is congruent to a complex curve in $\mathbb{C}^2$
without f\/lex points.

Section~\ref{s:examples} is devoted to the analysis of
explicit examples of applicable non-totally umbilic elliptic
Lagrangian immersions.

\newpage

\section{Preliminaries}\label{s:pre}

\subsection{Af\/f\/ine symplectic frames and structure equations}

Let $\mathbb{R}^4_\Omega$ denote $\mathbb{R}^4$ with the {\it standard symplectic form}
\[
  \Omega(X,Y)= {}^t\hskip-1pt X  J Y,\qquad  X,Y\in \mathbb{R}^4,
\]
where
\[
 J = (J_{ij}) =\left(
     \begin{array}{cc}
       0 & I_2 \\
        - I_2 & 0 \\
         \end{array}
          \right),
\]
$I_2$ being the $2\times 2$ identity matrix. The {\it linear
symplectic group} is
\[
  \mathrm{Sp}(4,\mathbb{R})=\{\mathbf{X} \in \mathrm{GL}(4,\mathbb{R}) : {}^t\hskip-1pt  \mathbf{X} J \mathbf{X} = J \}.
\]
An element $\mathbf{X} \in \mathrm{Sp}(4,\mathbb{R})$ is of the form
\[
 \mathbf{X} =\left(
    \begin{array}{cc}
     A & B \\
      C & D \\
       \end{array}
        \right),
\]
where $A$, $B$, $C$ and $D$ are $2\times 2$ matrices such that
\begin{gather*}
 {}^t\hskip-1pt  A  C-{}^t\hskip-1pt  C A =0, \qquad {}^t\hskip-1pt  B D-{}^t\hskip-1pt  D B = 0,
  \qquad {}^t\hskip-1pt  A D-{}^t\hskip-1pt  C B = I_{2}.
   \end{gather*}
The Lie algebra of $\mathrm{Sp}(4,\mathbb{R})$ is
\[
 \mathfrak{sp}(4,\mathbb{R})=\left\{
  \mathbf{x}(a,b,c)=\left(
    \begin{array}{cc}
       a & b \\
        c & -{}^t\hskip-1pt  a \\
         \end{array}
          \right) \, :
          \, a\in \mathfrak{gl}(2,\mathbb{R}),\, c,\, b\in \mathrm{S}(2,\mathbb{R}) \right\},
\]
where $\mathrm{S}(2,\mathbb{R})$ is the vector space of $2\times 2$ symmetric matrices
 endowed with the inner product of signature $(2,1)$ def\/ined by
\[
  (b,b) = -\det b, \qquad  b\in \mathrm{S}(2,\mathbb{R}).
\]
The {\it affine symplectic group} $\mathbb{R}^4 \rtimes \mathrm{Sp}(4,\mathbb{R})$
is represented in $\mathrm{GL}(5,\mathbb{R})$ by
\[
   \mathcal{S}(4,\mathbb{R})=
    \left\{\mathbf{S}(P,\mathbf{X})=\left(
     \begin{array}{cc}
       1 & 0 \\
         P & \mathbf{X} \\
         \end{array}
          \right) \,:\, P \in \mathbb{R}^4,\, \mathbf{X} \in \mathrm{Sp}(4,\mathbb{R}) \right\}.
\]
The Lie algebra of $\mathcal{S}(4,\mathbb{R})$ is given by
\[
  \mathfrak{s}(4,\mathbb{R})=\left\{\mathrm{S}(p,\mathbf{x})
     =\left(
       \begin{array}{cc}
         0 & 0 \\
          p & \mathbf{x} \\
           \end{array}
           \right)\, : \, p\in \mathbb{R}^4,\mathbf{x}\in \mathfrak{s}(4,\mathbb{R}) \right\}.
\]
If $\mathbf{S}(P,\mathbf{X})\in \mathcal{S}(4,\mathbb{R})$, we let $X_j$, $j=1,\dots,4$,
denote the column vectors of $\mathbf{X}$.

By an {\it affine symplectic frame} is meant a point $P\in \mathbb{R}^4$ and
four vectors $X_1, \dots, X_4 \in \mathbb{R}^4$ such that $\mathbf{X} =
(X_1,X_2,X_3,X_4) \in \mathrm{Sp}(4,\mathbb{R})$.
Upon choice of a reference frame,
the manifold of all af\/f\/ine symplectic frames of $\mathbb{R}^4_\Omega$
may be identif\/ied with the group $\mathcal{S}(4,\mathbb{R})$.

Consider the tautological projection maps
\begin{gather}
  \mathbf{P} : \ \mathcal{S}(4,\mathbb{R})\ni \mathbf{S}( P,\mathbf{X})\mapsto P\in \mathbb{R}^4, \label{ASS2}\\
   \mathbf{X}_j : \ \mathcal{S}(4,\mathbb{R})\ni \mathbf{S}( P,\mathbf{X}) \mapsto X_j\in \mathbb{R}^4,\nonumber
      \end{gather}
for $j=1,\dots,4$. Then
\[
  d\mathbf{P}=\tau^j \mathbf{X}_j, \qquad d\mathbf{X}_j
   =\theta_j^i \mathbf{X}_i,
\]
where the $\tau^j$, $\theta^i_j$ are the left-invariant Maurer--Cartan forms
of $\mathcal{S}(4,\mathbb{R})$. Exterior dif\/ferentiation of
\[
  \Omega(\mathbf{X}_i,\mathbf{X}_j)=J_{ij}
\]
yields that $\theta$ is $\mathfrak{sp}(4,\mathbb{R})$-valued. We then write
\[
 \theta =
  \left(\begin{array}{cc}
    \alpha & \beta \\
      \gamma & -{}^t\hskip-1pt \alpha \\
        \end{array}
         \right),
\]
where
\[
 \alpha =
   \left(\begin{array}{cc}
     \alpha^1_1 & \alpha^1_2 \\
       \alpha^2_1 & \alpha_2^2 \\
        \end{array}
           \right),
\qquad \beta =
 \left(\begin{array}{cc}
   \beta^1_1 & \beta^1_2 \\
     \beta^1_2 & \beta^2_2 \\
       \end{array}\right),
\qquad \gamma=
 \left(\begin{array}{cc}
    \gamma^1_1 & \gamma^1_2\\
      \gamma^1_2 & \gamma^2_2 \\
        \end{array}\right).
\]
Note that
\[
 \big(\tau^1,\tau^2,\tau^3,\tau^4,\alpha^1_1,\alpha^2_2,\alpha^1_2,
  \alpha^2_1,\gamma^1_1,\gamma^2_2,\gamma^1_2,\beta^1_1,\beta^2_2,\beta^1_2\big)
\]
is a basis for the vector space of the left-invariant 1-forms of
$\mathcal{S}(4,\mathbb{R})$. Moreover, for each $\mathbf{S}(0,\mathbf{Y}) \in \mathcal{S}(4,\mathbb{R})$, we
have
\begin{gather*}
 R^*_{\mathbf{S}(0,\mathbf{Y})}\left(\begin{array}{c}\tau^1\\\tau^2\\\tau^3\\\tau^4\end{array}\right)
   =\mathbf{Y}^{-1}\left(\begin{array}{c}\tau^1\\\tau^2\\\tau^3\\\tau^4\end{array}\right),
  \qquad R^*_{\mathbf{S}(0,\mathbf{Y})}(\theta)=\mathbf{Y}^{-1}\theta \mathbf{Y} + \mathbf{Y}^{-1}d\mathbf{Y}.
   \end{gather*}
Exterior dif\/ferentiation of
\[
  d\mathbf{P}=\tau^j\mathbf{X}_j,\qquad d\mathbf{X}_i=\theta^j_i\mathbf{X}_j
\]
yields the {\it structure equations} of the af\/f\/ine symplectic group:
\begin{gather} 
 \left\{\begin{array}{llll}
d\tau^1=-\alpha^1_1\wedge \tau^1-\alpha^1_2\wedge \tau^2-\beta^1_1\wedge \tau^3-\beta^1_2\wedge \tau^4,\\
d\tau^2=-\alpha^2_1\wedge \tau^1-\alpha^2_2\wedge \tau^2-\beta^1_2\wedge \tau^3-\beta^2_2\wedge \tau^4,\\
d\tau^3=-\gamma^1_1\wedge \tau^1-\gamma^1_2\wedge \tau^2+\alpha^1_1\wedge \tau^3+\alpha^2_1\wedge \tau^4,\\
d\tau^4=-\gamma^1_2\wedge \tau^1-\gamma^2_2\wedge
\tau^2+\alpha^1_2\wedge \tau^3+\alpha^2_2\wedge \tau^4,
\end{array}\right.\nonumber
\\
\label{MCA}
 \left\{\begin{array}{llll}
d\alpha^1_1=
 -\alpha^1_2\wedge\alpha^2_1-\beta^1_1\wedge \gamma^1_1-\beta^1_2\wedge \gamma^1_2,\\
d\alpha^2_2=
  -\alpha^2_1\wedge \alpha^1_2-\beta^1_2\wedge\gamma^1_2-\beta^2_2\wedge \gamma^2_2,\\
d\alpha^2_1=
  (\alpha^1_1-\alpha^2_2)\wedge \alpha^2_1-\beta^1_2\wedge \gamma^1_1-\beta^2_2\wedge \gamma^1_2,\\
d\alpha^1_2=
 (\alpha^2_2-\alpha_1^1)\wedge \alpha^1_2-\beta^1_1\wedge \gamma^1_2-\beta^1_2\wedge \gamma^2_2,
\end{array}\right.
\\
 \left\{\begin{array}{llll}
d\gamma^1_1=2\alpha^1_1\wedge \gamma^1_1+2\alpha^2_1\wedge \gamma^1_2,\\
d\gamma^2_2=2\alpha^2_2\wedge \gamma^2_2+2\alpha^1_2\wedge \gamma^1_2,\\
d\gamma^1_2=(\alpha^1_1+\alpha^2_2)\wedge
\gamma^1_2+\alpha^2_1\wedge \gamma^2_2+\alpha^1_2\wedge
\gamma^1_1,
\end{array}\right.\nonumber
\end{gather}
and
\begin{gather}\label{MCB}
 \left\{\begin{array}{llll}
d\beta^1_1=-2\alpha^1_1\wedge \beta^1_1-2\alpha^1_2\wedge \beta^1_2,\\
d\beta^2_2=-2\alpha^2_2\wedge \beta^2_2-2\alpha^2_1\wedge \beta^1_2,\\
d\beta^1_2=-(\alpha^1_1+\alpha^2_2)\wedge\beta^1_2-\alpha^2_1\wedge
\beta^1_1-\alpha^1_2\wedge \beta^2_2.
\end{array}\right.
\end{gather}

\subsection{Oriented Lagrangian planes}\label{ss:olp}

Let $\Lambda_2^+$ denote the set of oriented Lagrangian subspaces of
$\mathbb{R}^4$, i.e. the set of oriented 2-dimensional linear subspaces $V
\subset \mathbb{R}^4$ such that $\Omega_{| {V}}=0$. The space $\Lambda_2^+$
is a smooth manifold dif\/feomorphic to $S^2\times S^1$. The standard
action of $\mathrm{Sp}(4,\mathbb{R})$ on $\mathbb{R}^4$ induces an action on $\Lambda_2^+$
which is transitive.
The projection map
\begin{gather*}
  \pi_{\Lambda}: \ \mathbf{X}\in \mathrm{Sp}(4,\mathbb{R})\mapsto [X_1\wedge X_2]\in
  \Lambda_2^+
    \end{gather*}
makes $\mathrm{Sp}(4,\mathbb{R})$ into a principal f\/iber bundle
with structure group
\[
\mathrm{Sp}(4,\mathbb{R})_1 = \left\{ \mathbf{X}(A,b)=\left(
      \begin{array}{cc}
         A & Ab \\
           0 & {}^t\hskip-1pt  A^{-1} \\
            \end{array}
             \right) \, : \,  \det A > 0,\, b\in \mathrm{S}(2,\mathbb{R}) \right\}.
\]
From this, it follows that the forms
$(\gamma^1_1,\gamma_2^2,\gamma^2_1)$ span the semibasic forms of the
projection $\pi_{\Lambda}$. Moreover, the symmetric quadratic form
${g}=-\gamma^1_1\gamma^2_2+(\gamma^2_1)^2$ and the exterior 3-form
$\gamma^1_1\wedge \gamma^2_1\wedge \gamma^2_2$ are well def\/ined on
$\Lambda_2^+$, up to a positive multiple. They determine a conformal
structure of signature $(2,1)$ and an orientation, respectively.

\begin{remark}
The group $\mathrm{Sp}(4,\mathbb{R})$ is a covering group of the identity
component of the conformal Lorentz group $\mathrm{O}(3,2)$. This
fact implies that $\Lambda_2^+$ can be identif\/ied with the quotient
of the nullcone of $\mathbb{R}^{3,2}$ by the action given by positive scalar
multiplications. Therefore, $\Lambda_2^+$ can be seen as the {\it
conformal compactification} of oriented, time-oriented af\/f\/ine
Minkowski space $\mathbb{R}^{2,1}\cong \mathrm{S}(2,\mathbb{R})$ (cf.~\cite{BCDGM,GSkepler}).
Any element of $\Lambda_2^+$ can be represented by a $4\times 2$
matrix
\[
X = \left(\begin{array}{c} A_1\\ A_2 \end{array} \right)
\]
of rank 2 such that ${}^t\hskip-1pt  A_1 A_2 = {}^t\hskip-1pt  A_2 A_1$.
Let $U_0 \subset \Lambda_2^+$ be the set of
all Lagrangian subspaces $[X_1\wedge X_2]$ associated to ordered
pairs of vectors
\[
 X_1= {}^t\hskip-1pt  \big(x_1^1,\dots,x_1^4\big), \qquad
  X_2= {}^t\hskip-1pt \big(x_2^1,\dots,x_2^4\big)
\]
such that
\[
 \det\left(\begin{array}{cc}
     x^1_1 & x_2^1 \\
      x^2_1 & x^2_2 \\
       \end{array}
        \right)> 0.
\]
The set $U_0$ is a dense open subset of $\Lambda_2^+$. The map $S :
U_0 \to \mathrm{S}(2,\mathbb{R})\cong \mathbb{R}^{2,1}$ given by
\[
 \mathrm{S}([X_1\wedge X_2])=
   \left(
     \begin{array}{cc}
         x^3_1 & x_2^3 \\
          x^4_1 & x^4_2 \\
           \end{array}
            \right)
  \left(
   \begin{array}{cc}
     x^1_1 & x_2^1 \\
      x^2_1 & x^2_2 \\
       \end{array}
         \right)^{-1}
\]
is a conformal dif\/feomorphism, and $(U_0,S)$ is a local chart of
$\Lambda_2^+$.
Note that
\[
 U_0 = \left\{ V \in \Lambda_2^+ \, : \; V \cap W_0 = \{0\}\right\},
\]
the set of Lagrangian subspaces which are transversal to $W_0 =
\mathrm{span}\,\{e_3,e_4\}$.\footnote{Here $e_1$, $e_2$, $e_3$, $e_4$ denote the
standard basis of $\mathbb{R}^4$.} Let
\[
 U_1 = \left\{ V \in \Lambda_2^+ \, : \; \dim (V \cap W_0) =
  1\right\}.
\]
Then $U_1 \cup \{W_0\}$ is a closed subset of $\Lambda_2^+$, which
can be interpreted as the ideal boundary of $\mathbb{R}^{2,1}\cong U_0$.
Note that $U_0$ and $U_1$ are orbits of the closed subgroup of
$\mathrm{Sp}(4,\mathbb{R})$ that preserves the subspace $W_0$.
 \end{remark}

\subsection{Lagrangian immersions}

Let $f : M \to \mathbb{R}^4$ be a smooth immersion of a connected
2-manifold, oriented by a volume element~$\mathcal{A}$. For any
$q\in M$, let $\mathcal{T}_f(q)$ be the 2-plane $df(T_q M)$
translated to the origin and equipped with the orientation induced
by $- df(\mathcal{A}|_q)$.\footnote{One reason for choosing the
opposite orientation is that totally umbilic elliptic Lagrangian
immersions are then described in terms of holomorphic functions
instead of anti-holomorphic ones (cf.\ Sections~\ref{section4.4} and~\ref{section4.5}).}

\begin{definition}
The immersion $f : M \to \mathbb{R}^4$ is called \textit{Lagrangian} if
$\mathcal{T}_f(q) \in \Lambda_2^+$, for each $q\in M$. The resulting
map
\[
  \mathcal{T}_f : q \in M\mapsto \mathcal{T}_f(q) \in \Lambda_2^+
\]
is the (\textit{symplectic}) \textit{Gauss map} of the
Lagrangian immersion $f$.
\end{definition}

\begin{definition}
A Lagrangian immersion is said of \textit{general type}
(respectively, \textit{special type}) if its  Gauss map is an
immersion (respectively, has constant rank one). Moreover, a
Lagrangian immersion of general type is said {\it elliptic}
(respectively, {\it hyperbolic}, {\it parabolic}) if its Gauss map
is a spacelike (respectively, timelike, lightlike) immersion into~$\Lambda_2^+$.
\end{definition}

\section{Moving frames on elliptic Lagrangian immersions}\label{s:reduction}

\subsection{Frame reductions and dif\/ferential invariants}

In this section, $f : M\to \mathbb{R}^4$ will denote an elliptic Lagrangian
immersion with Gauss map $\mathcal{T}_f$. Since we are working with
immersions which are not necessarily one-to-one, it is not
restrictive to assume that $M$ is simply connected.

\begin{definition}
A {\it symplectic frame field along} $f$ is a smooth map
\[
 \mathbf{S} = \mathbf{S}(f,\mathbf{X}) : \ U \to \mathcal{S}(4,\mathbb{R})
   \]
from an open connected subset $U\subset M$ such that the projection
\eqref{ASS2} composed with $\mathbf{S}$ is $f$, i.e.\ $\mathbf{P} \circ
\mathbf{S}(f,\mathbf{X}) = \mathbf{S}(f,\mathbf{X})\cdot O = f$.
\end{definition}

Following the usual practice in the method of moving frames we will
write $\phi$ instead of $\mathbf{S}^*(\phi)$ to denote forms on $U$ which
are pulled back from $\mathcal{S}(4,\mathbb{R})$ by the moving frame $\mathbf{S}$. A
symplectic frame f\/ield along $f$ gives a local representation of the
derivative map,
\[
  df =\sum_{j=1}^{4}\tau^j\mathrm{X}_j.
   \]
Let $\mathbf{S}:U\to \mathcal{S}(4,\mathbb{R})$ be a symplectic frame f\/ield along $f$. Any other
symplectic frame f\/ield along~$f$ on~$U$ is given by
\[
  \widetilde{\mathbf{S}}=\mathbf{S} \mathbf{Y},
  \]
where $\mathbf{Y}:U\to \mathrm{Sp}(4,\mathbb{R})$ is a smooth map\footnote{Here we are
identifying an element $\mathbf{S}(0,\mathbf{Y}) \in \mathcal{S}(4,\mathbb{R})$ with
$\mathbf{Y}\in\mathrm{Sp}(4,\mathbb{R})$, and hence an element $\mathrm{S}(0,\mathbf{x})
\in \mathfrak{s}(4,\mathbb{R})$ with $\mathbf{x} \in \mathfrak{sp}(4,\mathbb{R})$.}.
If $\Theta$ and $\widetilde{\Theta}$ are the pull-backs of the Maurer--Cartan
form of $\mathcal{S}(4,\mathbb{R})$ by $\mathbf{S}$ and $\widetilde{\mathbf{S}}$, respectively, then
\begin{gather}\label{0F1}
  \widetilde{\Theta}=\mathbf{Y}^{-1}\Theta \mathbf{Y}+\mathbf{Y}^{-1}d\mathbf{Y}.
   \end{gather}

\subsubsection{First order frame f\/ields}

A symplectic frame f\/ield $\mathbf{S}=\mathbf{S}(f,\mathbf{X}) : U\to \mathcal{S}(4,\mathbb{R})$ along $f$
is of {\it first order} if
\[
 \mathcal{T}_f(q) = [\mathrm{X}_1|_q\wedge \mathrm{X}_2|_q], \qquad \text{for each $q\in U$}.
\]
It is clear that f\/irst order frame f\/ields exist locally.
With respect to a f\/irst order frame, the 1-forms
$\tau^3$ and $\tau^4$ vanish identically on $U$ and
\[
  df=\tau^1\mathrm{X}_1+\tau^2\mathrm{X}_2,
   \]
where $\tau^1\wedge \tau^2<0$. If $\mathbf{S}$ is a f\/irst order frame f\/ield
along $f$ on $U$, then any other is given by $\widetilde{\mathbf{S}}=\mathbf{S}
\mathbf{Y}(A,b)$, where $\mathbf{Y}(A,b) : U \to \mathrm{Sp}(4,\mathbb{R})_1$ is a
smooth map into the subgroup $\mathrm{Sp}(4,\mathbb{R})_1$ introduced in Section
\ref{ss:olp}.

Calculated with respect to a f\/irst order frame f\/ield, the quadratic
form $g = -\gamma^1_1\gamma^2_2+(\gamma^2_1)^2$ is positive
def\/inite. Under a change $\widetilde{\mathbf{S}}=\mathbf{S}\mathbf{Y}(A,b)$ of
f\/irst order frames, it transforms by
\[
  -\widetilde{\gamma}^1_1\widetilde{\gamma}^2_2+(\widetilde{\gamma}^2_1)^2
       =\det(A)^2\big(-\gamma^1_1\gamma^2_2+(\gamma^2_1)^2\big),
\]
and hence def\/ines a conformal structure on $M$. On $M$, we will consider
the unique complex structure compatible with the given orientation and the
conformal structure def\/ined by $g$. The 1-forms
\[
  \omega^1:=\tfrac{1}{2}\big(\gamma^1_1-\gamma^2_2\big),\qquad \omega^2:=\gamma^2_1
    \]
are everywhere linearly independent and
\[
  \omega=\omega^1+i\omega^2
   \]
is complex-valued of bidegree $(1,0)$. Moreover, there exist smooth
functions $\ell_1,\ell_2:U\to \mathbb{R}$ such that
\[
 \tfrac{1}{2}\big(\gamma^1_1+\gamma^2_2\big)=\ell_1\omega^1+\ell_2\omega^2,\qquad
  \ell_1^2+\ell_2^2 < 1.
    \]

\subsubsection{Second order frame f\/ields}
A f\/irst order frame f\/ield along $f$ is of {\it second order} if it satisf\/ies
\[
  \tfrac{1}{2}\big(\gamma^1_1+\gamma^2_2\big)=0.
\]

\begin{lemma}
About any point of $M$ there exists a second order frame field along
$f$.
 \end{lemma}

\begin{proof}
Let $\mathbf{S}:U\to \mathcal{S}(4,\mathbb{R})$ be any f\/irst order frame f\/ield along $f$.
Let $\phi:U\to (-\pi/2,\pi/2)$ be the smooth function def\/ined by
\[
  \sin(\phi)=-\frac{\ell_2}{\sqrt{1-\ell_1^2}}
\]
and let $A:U\to \mathrm{GL}_+(2,\mathbb{R})$ be given by
\[
  A=\left(
     \begin{array}{cc}
       \sqrt{\frac{1-\ell_1}{2}}\cos(\phi) & \sqrt{\frac{1-\ell_1}{2}}\sin(\phi) \\
         0 &  \sqrt{\frac{1+\ell_1}{2}}
         \end{array}
          \right).
           \]
Then $\mathbf{Y}(A,0):U\to \mathrm{Sp}(4,\mathbb{R})_1$ is a smooth map. Let
$\widetilde{\mathbf{S}}=\mathbf{S} \mathbf{Y}(A,0)$. The transformation rule~(\ref{0F1}) gives $\widetilde \gamma = {}^t\hskip-1pt  A \gamma A$, from which we
compute $\widetilde \gamma^1_1 + \widetilde \gamma^2_2 = 0$. This
means that $\widetilde{\mathbf{S}}$ is a second order frame f\/ield along
$f$.
\end{proof}

If $\mathbf{S}$ and $\widetilde{\mathbf{S}}$ are second order frame f\/ields on $U$,
then $\widetilde{\mathbf{S}} = \mathbf{S} \mathbf{Y}(A,b)$, where $A$ is a
$2\times 2$ matrix such that ${}^t\hskip-1pt  A  A=r^2I_2$ and $\mathrm{det}(A) >
0$. The group of all such $A$ identif\/ies with the multiplicative
group $\mathbb{C}^\ast$ of nonzero complex numbers by
\[
  {\mathbb{C}}^\ast \ni re^{is}
   \mapsto r\left(
      \begin{array}{cc}
        \cos(s) & -\sin(s) \\
           \sin(s) & \cos(s) \\
             \end{array}
               \right).
                  \]
Let $\mathbf{S}$ be a second order frame f\/ield along $f$. We know that the
complex-valued 1-form $\omega=\omega^1+i\omega^2$
is of bidegree $(1,0)$ and never zero.

Dif\/ferentiating the equations $\tau^3=\tau^4=0$ and using the
structure equations, we f\/ind that
\[
 \left\{\begin{array}{llll}
  \tau^1\wedge \omega^1+\tau^2\wedge \omega^2=0,\\
   \tau^1\wedge \omega^2-\tau^2\wedge \omega^1=0,\\
    \end{array}\right.
      \]
which, by Cartan's Lemma, implies that
\[
  \tau :=\tau^1-i\tau^2 = t\omega,
\]
for some smooth function $t = t_1 +it_2 : U \to \mathbb{C}$.
The 1-form $\tau$ is of bidegree $(1,0)$ and never zero.

Under a change of second order frame f\/ields $\widetilde{\mathbf{S}}=\mathbf{S}
\mathbf{Y}(re^{is},b)$, the forms $\omega$ and $\tau$ transform~by
\begin{gather}\label{2F1}
  \widetilde{\omega}=r^2e^{-2is}\omega,\qquad \widetilde{\tau}=r^{-1}e^{is}\tau.
   \end{gather}
This implies that
\[
 \widetilde{\tau}^2\widetilde{\omega} = \tau^2 \omega.
  \]
Thus $\tau^2 \omega$ is independent of the
choice of second order frame f\/ield and is never zero.

\begin{definition}
The cubic dif\/ferential form $\mathcal{F} : = \tau^2 \omega$,
globally def\/ined on $M$, is called
the \textit{Fubini cubic form} of the
Lagrangian immersion~$f$.
\end{definition}

 Next, consider
the complex-valued 1-form
\[
  \eta:=\eta_1-i\eta_2=\tfrac{1}{2}\big(\alpha^1_1-\alpha^2_2\big)
    -\tfrac{i}{2}\big(\alpha^2_1+\alpha^1_2\big).
     \]
Dif\/ferentiation of $\gamma^1_1+\gamma^2_2=0$, combined with
the structure equations, gives
\[
  \eta_1\wedge \omega^1+\eta_2\wedge \omega^2=0,
    \]
which, by Cartan's Lemma, implies that
\[
  \eta = h\omega+\ell \bar{\omega}
    \]
for some smooth functions $h = h_1 +ih_2 : U \to \mathbb{C}$ and $\ell : U \to \mathbb{R}$.

\subsubsection{Third order frame f\/ields}

Under a change $\widetilde{\mathbf{S}}=\mathbf{S} \mathbf{Y}(re^{is},b)$ of second
order frame f\/ields, the 1-form $\eta$ transforms by
\[
   \widetilde{\eta}=e^{2is}\left(\eta -
   \frac{r^2}{2}\,
    \mathrm{tr}\,(b)\,\bar{\omega}\right).
     \]
This implies that locally there exist second order frame f\/ields satisfying
$\ell=0$. Such frame f\/ields are said of {\it third order}.

If $\widetilde{\mathbf{S}}$ and $\mathbf{S}$ are third order frame
f\/ields on $U$, then $\widetilde{\mathbf{S}} = \mathbf{S} \mathbf{Y}(re^{is},b)$
where $\mathrm{tr}\,(b)=0$. With respect to a third order frame f\/ield, the
1-form $\eta$ is of bidegree $(1,0)$ and transforms by
\begin{gather}\label{2F11}
   \widetilde{\eta}=e^{2is}\eta.
     \end{gather}

\begin{definition}
According to (\ref{2F1}) and (\ref{2F11}), the dif\/ferential forms
\[
  \mathcal{B}  =\bar{\eta}\eta,\qquad
    \mathcal{H}  =|t|^{2/3}\eta\omega,
       \]
are globally def\/ined on $M$, independent of the choice of third order frames.
We call $\mathcal{B}$ the \textit{Thomsen quadratic form} and
$\mathcal{H}$ the \textit{normalized Hopf differential} of $f$.
\end{definition}

\begin{remark}
The invariant form $\mathcal{B}$ amounts to the metric induced on
$M$ by the conformal Gauss map \cite{BrDG} (or central sphere
congruence \cite{Bl2}) of the spacelike conformal immersion
$\mathcal{T}_f : M\to \Lambda_2^+$. The invariant form $\mathcal{H}$
is instead a normalization of the Hopf dif\/ferential of
$\mathcal{T}_f$.
\end{remark}

\subsection{Adapted frame f\/ields}

Being simply connected, $M$ is
either the Riemann sphere $S^2$, the complex plane $\mathbb{C}$,
or the unit disk $\Delta$. By a \textit{complex parameter} on $M$ is meant a
complex coordinate chart $(U,z)$ in $M$ def\/ined on a maximal
contractible open subset $U$ of $M$: if $M=\mathbb{C}$, or $\Delta$, then
$U=M$; if $M=S^2$, then $U=S^2\setminus\{q\}$.

\begin{definition}
Let $(U,z)$ be a complex parameter on $M$.
A third order frame f\/ield $\mathbf{S} : U \to \mathcal{S}(4,\mathbb{R})$ along $f$
is \textit{adapted} to $(U,z)$ if
\begin{gather}\label{adapted-frame}
  \omega =\omega^1+i\omega^2=dz,\qquad \alpha^1_1
    +\alpha^2_2=\alpha^2_1-\alpha^1_2=0.
     \end{gather}
\end{definition}

\begin{proposition}
Let $(U,z)$ be a complex parameter on $M$. There exist
third order frame fields
along $f$ adapted to $(U,z)$.
If $\mathbf{S}$ and $\mathbf{S}'$ are two such frame fields, then $\mathbf{S}'=\pm \mathbf{S}$.
\end{proposition}

\begin{proof}
Let $q_0 \in U$ and let $\mathbf{S}' : V\to \mathcal{S}(4,\mathbb{R})$ be any third order
frame f\/ield def\/ined on a simply connected open neighborhood $V
\subset U$ of $q_0$. On $V$ there exist real-valued functions
$r,\vartheta:V\to \mathbb{R}$ such that $\omega'=r^2e^{-2i\vartheta}dz$.
Def\/ine $\mathbf{S}''=\mathbf{S}' \mathbf{Y}(r^{-2}e^{2i\vartheta},0)$. Then,
$\mathbf{S}''$ is a third order frame f\/ield such that $\omega''=dz$. Two
such frame f\/ields are related by
\[
  \widetilde{\mathbf{S}}''=\mathbf{S}'' \mathbf{Y}(\pm I_2,b),
 \]
where $b:V\to \mathrm{S}(2,\mathbb{R})$ is a smooth map such that $\mathrm{tr}\,(b)=0$.
From the structure equations we f\/ind that
\[
 \left\{\begin{array}{llll}
  \big(\alpha^1_1+\alpha^2_2\big)\wedge dx+\big(\alpha^2_1-\alpha^1_2\big)\wedge dy =0,\\
   \big(\alpha^1_1+\alpha^2_2\big)\wedge dy-\big(\alpha^2_1-\alpha^1_2\big)\wedge dx =0,
    \end{array}\right.
\]
which implies
\[
  (\alpha^1_1+\alpha^2_2)-i(\alpha^2_1-\alpha^1_2)=w\,dz,
\]
for some smooth function $w : V \to \mathbb{C}$.
If $\widetilde{\mathbf{S}}''=\mathbf{S}'' \mathbf{Y}(\pm I_2,b)$ then
\[
  \big(\tilde{\alpha}^1_1+\tilde{\alpha}^2_2\big)
    -i\big(\tilde{\alpha}^2_1-\tilde{\alpha}^1_2\big)=\big(\alpha^1_1+\alpha^2_2\big)
      -i\big(\alpha^2_1-\alpha^1_2\big)+ 2\big(b^1_1-ib^2_1\big)dz.
\]
If we choose $b$ such that $2(b^1_1-ib^2_1)=-w$, then
$\widetilde{\mathbf{S}}''$ satisf\/ies
\[
 \omega''=dz,\qquad \big(\tilde{\alpha}^1_1+\tilde{\alpha}^2_2\big)
  -i\big(\tilde{\alpha}^2_1-\tilde{\alpha}^1_2\big)=0.
\]
This shows that adapted frame f\/ields do exist locally near any point
of $U$. Moreover, two such frames $\mathbf{S}$ and $\mathbf{S}'$ on $V$ are
related by $\mathbf{S}'=\pm \mathbf{S}$. The existence of an adapted frame on $U$
follows from the fact that $U$ is simply connected.
\end{proof}

\subsubsection{Invariant functions and integrability equations}

Let $\mathbf{S}:U\to \mathcal{S}(4,\mathbb{R})$ be a frame f\/ield adapted to
 $(U,z)$. Then
\begin{gather}\label{t&h}
  \tau =t \,dz, \qquad \eta = h\,dz,
   \end{gather}
for smooth functions $t=t_1 +it_2: U\to \mathbb{C}$ and $h=h_1 +ih_2: U\to \mathbb{C}$.
Def\/ine $\upsilon$ and $\rho$ by
\[
 \upsilon:=\tfrac{1}{2}\big(\beta^1_1+\beta^2_2\big),\qquad \rho:=
   \tfrac{1}{2}\big(\beta^1_1-\beta^2_2\big)-i\beta^2_1.
\]
The structure equations (\ref{MCA}) imply
\begin{gather}\label{upsilon}
 \upsilon= h_{\bar{z}}dz+\bar{h}_zd\bar{z},
  \end{gather}
and
\begin{gather}\label{p}
 \rho = pdz + \big(h_1^2+h_2^2\big)d\bar{z} = pdz+ |h|^2 d\bar{z},
   \end{gather}
for some smooth function $p = p_1 + ip_2 : U\to \mathbb{C}$.

\begin{definition}
The functions $t$, $h$ and $p$ def\/ined by \eqref{t&h} and \eqref{p}
are called the {\it invariant functions} of the frame f\/ield along
$f$ adapted to $(U,z)$. These functions are subject to the equations
\begin{gather}\label{INTEQ}
\left\{\aligned
  &t_{\bar{z}} =\bar{t}\, h,\\
    &p_{\bar{z}} =2h\bar{h}_{z}-i |h|^2_z,\\
       & \bar p h -p\bar h = h_{\bar{z}\bar{z}} - \bar{h}_{zz}.
\endaligned
 \right.
         \end{gather}
\end{definition}

\begin{remark}
From \eqref{INTEQ}
it follows that
\begin{gather}\label{2F12}
 \bar{\partial}(\mathcal{F}) = 2\mathcal{N},
  \end{gather}
where $\mathcal{N}: = \bar{\tau}\tau \eta \omega$ is an invariant
form. The form $\mathcal{B}$, $\mathcal{H}$ and $\mathcal{N}$ are
completely determined by the Fubini cubic form.
\end{remark}

\begin{theorem}[Existence]\label{thm:fund}
Let $(U,z)$ be a complex parameter on $M$. Let $t : U \to \mathbb{C}$,
$h : U\to \mathbb{C}$ and $p : U \to \mathbb{C}$ be smooth functions satisfying
the equations $(\ref{INTEQ})$.
Then there exist an elliptic Lagrangian immersion
$f : U \to \mathbb{R}^4$, unique up to
affine symplectic transformation, and a~unique frame field $\mathbf{S}:U\to
\mathcal{S}(4,\mathbb{R})$ along $f$ adapted to $(U,z)$ whose invariant functions are~$t$,~$h$ and~$p$,~i.e.
\begin{gather*}
 \tau = tdz,\qquad \eta = hdz,\qquad \rho^{(1,0)} = pdz.
\end{gather*}
\end{theorem}

\begin{proof}
Let
\begin{gather*}
 \left(
   \begin{array}{c}
    \tau^1 \\
     \tau^2 \\
      \end{array}
       \right)=\left(
     \begin{array}{c}
     t_1dx-t_2dy \\
      -t_2dx-t_1dy \\
       \end{array}
        \right),\qquad \gamma=\left(
  \begin{array}{cc}
  dx & dy \\
   dy & -dx \\
   \end{array}
    \right),
\\
 \alpha=\left(
 \begin{array}{cc}
   h_1dx-h_2dy & -h_2dx-h_1dy \\
   -h_2dx-h_1dy &  -h_1dx+h_2dy\\
    \end{array}
     \right),\qquad
   \beta=\left(
    \begin{array}{cc}
     \beta^1_1 & \beta^2_1 \\
      \beta^2_1 & \beta^2_2 \\
       \end{array}
        \right),
 \end{gather*}
where
\[
 \tfrac{1}{2}\big(\beta^1_1+\beta^2_2\big)= h_{\bar{z}}dz + \bar{h}_zd\bar{z}
\]
and
\[
 \tfrac{1}{2}\big(\beta^1_1-\beta^2_2\big)-i\beta^2_1 = pdz+ |h|^2 d\bar{z}.
\]
Then
\begin{gather*}
 \Theta=\left(
  \begin{array}{ccc}
    0 & 0 & 0 \\
    \left(\begin{array}{c} \tau^1\\ \tau^2\end{array}\right) & \alpha & \beta \\
     0 & \gamma & -{}^t\hskip-1pt  \alpha \\
     \end{array}
      \right)
       \end{gather*}
is a smooth 1-form on $U$ with values in $\mathfrak{s}(4,\mathbb{R})$.
By (\ref{INTEQ}),
\[
  d\Theta= -\Theta \wedge \Theta.
    \]
Therefore, by the Cartan--Darboux theorem, there exists a smooth map
$\mathbf{S} : U \to \mathcal{S}(4,\mathbb{R})$ such that $\mathbf{S}^{-1}d\mathbf{S} = \Theta$, unique up
to left multiplication by an element of $\mathcal{S}(4,\mathbb{R})$. The elliptic
Lag\-ran\-gian immersion $f:U\to \mathbb{R}^4$ def\/ined by $f=\mathbf{S}\cdot O$ has the
required properties.
\end{proof}

\begin{remark}[Uniqueness]
If $f$, $\widetilde{f} : U \to \mathbb{R}^4$ are two elliptic Lagrangian immersions
inducing the same invariant functions, then, by the Cartan--Darboux uniqueness
theorem, there exists a~symplectic motion $\mathbf{S} : \mathbb{R}^4 \to \mathbb{R}^4$ such that
$\widetilde{f} = \mathbf{S} f$.
\end{remark}

\section{Symplectic applicability}\label{s:applicability}

In this section we investigate the extent to which the invariants are
really needed to determine an elliptic Lagrangian immersion
up to symmetry. We take the point of view of
surface applicability as in classical projective dif\/ferential
geometry and Lie sphere geometry \cite{Fu,Ca1,Ca2,Bl2,MNtohoku,JM}.

\subsection{Generic Lagrangian immersions}

Suppose that $\mathcal{H}$, the normalized Hopf dif\/ferential, is
never zero. Then the Hermitian form
\[
\mathcal{P}=\left((\bar{h}_{z}{\bar{h}}^{-1})_{\bar{z}} -
 ({h}_{\bar z}{{h}}^{-1})_{{z}}\right)dz d\bar{z}
\]
is well def\/ined, independent of the choice of the complex parameter
$(U,z)$ and of the adapted frame (the reason for considering this
form will be clear in the proof of the next proposition).

\begin{definition}
An elliptic Lagrangian immersion $f : M \to \mathbb{R}^4$ is called
{\it generic} if $\mathcal{H}(q)\neq 0$ and $\mathcal{P}(q)\neq 0$,
for each $q \in M$.
\end{definition}

\begin{proposition}\label{prop:gen}
A generic elliptic Lagrangian immersion is uniquely determined,
up to affine symplectic transformation, by the
Fubini cubic form.
\end{proposition}

\begin{proof}
Fix a complex parameter $(U,z)$ and choose an adapted
frame f\/ield $\mathbf{S}:U\to \mathcal{S}(4,\mathbb{R})$. If $h\neq 0$, then (\ref{INTEQ})
implies the existence of a unique real-valued function
$s:U\to \mathbb{R}$ such that
\[
 p=hs + \mathcal{D}_2(h),
  \]
where
\[
 \mathcal{D}_2(h)= \frac{1}{2\bar{h}}\left(\bar{h}_{zz} -h_{\bar{z}\bar{z}}\right).
\]
Dif\/ferentiating this identity and using
\[
  p_{\bar{z}}=2h\bar{h}_{z}-i |h|^2_z,
    \]
we get
\begin{gather}\label{CHE}
 ds=-s\left(\frac{1}{h} h_{\bar{z}}d\bar{z}
  +\frac{1}{\bar{h}}\bar{h}_{z}dz\right)
   -\left(\mathcal{D}_3(h)d\bar{z}
    +\overline{\mathcal{D}_3(h)}dz\right),
       \end{gather}
where
\[
 \mathcal{D}_3(h) =
 \frac{1}{h}\left[(\mathcal{D}_2(h))_{\bar{z}}-2h\bar{h}_{z}
  +i |h|^2_z\right].
   \]
Dif\/ferentiation of (\ref{CHE}) gives
\begin{gather*}
 \mathcal{P}_2(h)s+\mathcal{D}_4(h)=0,
  \end{gather*}
where
\[
 \mathcal{P}_2(h)= (\bar{h}_{z}{\bar{h}}^{-1})_{\bar{z}} -
 ({h}_{\bar z}{{h}}^{-1})_{{z}}
     \]
and
\[
 \mathcal{D}_4(h)= \overline{\mathcal{D}_3(h)}_{\bar{z}}
   -\mathcal{D}_3(h)_z -
     {\bar{h}}^{-1}\bar{h}_z \mathcal{D}_3(h)
      +{h}^{-1} h_{\bar{z}}\overline{\mathcal{D}_3(h)}.
        \]
If $\mathcal{P}_2(h)\neq 0$, then
\[
  s=-\frac{\mathcal{D}_4(h)}{\mathcal{P}_{2} (h)}.
   \]
Therefore, for a generic elliptic Lagrangian immersion (i.e.\ $h\neq
0$ and $\mathcal{P}_2(h)\neq 0$), the functions~$h$ and~$p$ are
determined by the Fubini cubic form. Thus, a generic $f$ is uniquely
determined, up to symplectic congruence, by its Fubini cubic form.
\end{proof}

\subsection{Applicable elliptic Lagrangian immersions}

In this section we investigate the special class of surfaces which are not
determined by the Fubini cubic form alone.

\begin{definition}
Two noncongruent elliptic Lagrangian immersions $f$, $\widetilde{f}
: M\to \mathbb{R}^4$ are \textit{applicable} if they induce the same Fubini
cubic forms, $\mathcal{F}= \widetilde{\mathcal{F}}$. An elliptic
Lagrangian immersion $f : M \to \mathbb{R}^4$ is \textit{applicable} if
there exists an elliptic Lagrangian immersion $\widetilde{f} : M\to
\mathbb{R}^4$ such that $f$ and $\widetilde{f}$ are applicable. Moreover, an
elliptic Lagrangian immersion $f$ is \textit{rigid} if any elliptic
Lagrangian immersion inducing the same Fubini cubic form as that of
$f$ is congruent to $f$.
\end{definition}

\begin{remark}
According to Proposition \ref{prop:gen}, any generic elliptic Lagrangian
immersion is rigid.
\end{remark}

\begin{theorem}\label{thm:applicable}
An elliptic Lagrangian immersion $f:M\to \mathbb{R}^4$ is applicable if and
only if there exists a non-zero holomorphic quadratic differential
$\mathcal{Q}$ on $M$ such that $\mathcal{H}\wedge_{\mathbb{R}}
\mathcal{Q}=0$.\footnote{$\mathcal{H}\wedge_{\mathbb{R}}
\mathcal{Q}=0$ means that $\mathcal{H}$ and $\mathcal{Q}$ are
linearly dependent over the reals.}
\end{theorem}

\begin{proof}
First, we consider the case
$\mathcal{H}\equiv 0$. Then, $\mathcal{F}$ is a never vanishing
holomorphic cubic dif\/ferential. Since $M$ is simply connected, $M$
is either the complex plane or the unit disk. In both cases, there
exist a global complex coordinate $z : M \to \mathbb{C}$ and a global
adapted frame f\/ield $\mathbf{S} : M \to \mathcal{S}(4,\mathbb{R})$ with invariant functions
$t : M \to \mathbb{C}$, $h = 0$ and $p : M\to \mathbb{C}$. Consider a non-zero
holomorphic dif\/ferential $\mathcal{Q}=\lambda \,dz^2$ and set
$\widetilde{p}=p-\lambda$. The functions $t$, $h= 0$ and
$\widetilde{p}$ satisfy equations (\ref{INTEQ}). The corresponding
$\widetilde{f}$ has $\widetilde{\mathcal{F}} =\mathcal{F}$, but is
not congruent to $f$.

We now examine the case $\mathcal{H}\neq 0$.
Let $f$, $\widetilde{f}$ be two
noncongruent elliptic Lagrangian immersions inducing the same Fubini
cubic forms.
Let $\mathbf{S}$, $\widetilde{\mathbf{S}}:U\to \mathcal{S}(4,\mathbb{R})$ be the corresponding adapted
frame f\/ields with respect to a f\/ixed complex parameter $(U,z)$.
Since $t = \widetilde{t}$, equations~(\ref{INTEQ}) imply that
$h = \widetilde{h}$ and that $p-\widetilde{p}$ is a non-zero
holomorphic function satisfying
\[
  (p_1-\widetilde{p}_1)h_2-(p_2-\widetilde{p}_2)h_1=0.
    \]
Then, $\mathcal{Q}_U = (p{-}\widetilde{p})dz^2$ is a non-zero
holomorphic quadratic dif\/ferential on~$U$ such that
$\mathcal{H}\wedge_{\mathbb{R}} \mathcal{Q}_U= 0$. Let
$(\widehat{U},\widehat{z})$ be another complex parameter and let
${\widehat{\mathcal{Q}}}_{\hat{U}}$ be the corresponding holomorphic
quadratic dif\/ferential. The equations
$\mathcal{H}\wedge_{\mathbb{R}}\mathcal{Q}_U=0$ and
$\mathcal{H}\wedge_{\mathbb{R}}\widehat{\mathcal{Q}}_{\hat U}=0$ imply that
$\mathcal{Q}=\lambda_{U\hat{U}} \widehat{Q}$ on $U\cap
\widehat{U}\neq \varnothing$, where $\lambda_{U\hat{U}}$ is a
real-valued, locally constant function. The functions
$\{\exp{\lambda_{U\hat{U}}}\}$, def\/ined on each ordered pair
$(U,\widehat{U})$ when $U\cap \widehat{U}\neq \varnothing$, generate
an $\mathbb{R}^+$-valued 1-cocycle on~$M$, which is trivial because $M$ is
simply connected. Therefore, there exists a globally def\/ined
non-zero holomorphic quadratic dif\/ferential $\mathcal{Q}'$ such that
$\mathcal{Q}_U=\mathcal{Q}'|_U$, for every $(U,z)$.

Conversely, suppose there exists a non-zero holomorphic quadratic
dif\/ferential $\mathcal{Q}$ such that $\mathcal{H}\wedge_{\mathbb{R}}
\mathcal{Q}=0$. This implies that $M$ is either $\mathbb{C}$ or the unit
disk. Let $z$ be a global complex coordinate on $M$ and let
$\mathbf{S}:M\to \mathcal{S}(4,\mathbb{R})$ be the associated adapted frame f\/ield. Set
$\mathcal{Q}=w\,dz^2$, where $w=w_1 +iw_2$ is a holomorphic non-zero
function. Now, the invariant functions~$t$,~$h$ and~$p$, calculated
with respect to~$\mathbf{S}$, satisfy (\ref{INTEQ}). Since
$\mathcal{H}\wedge_{\mathbb{R}} \mathcal{Q}=0$, that is $w_1h_2 -w_2h_1 =0$,
also $t$, $h$ and $\widetilde{p} = p-w$ satisfy~(\ref{INTEQ}). The
corresponding elliptic Lagrangian immersions $f$,
$\widetilde{f}:M\to \mathbb{R}^4$ are applicable.
\end{proof}

\begin{remark}
If $\mathcal{H}\neq 0$, then equation
$\mathcal{H}\wedge_{\mathbb{R}}\mathcal{Q}=0$ determines $\mathcal{Q}$ up to
a real constant multiple. For each $\lambda\in \mathbb{R}$, def\/ine the
invariant functions $t_\lambda = t$, $h_\lambda = h$ and $p_\lambda
= p -\lambda$. It is clear that $t_\lambda$, $h_\lambda$ and
$p_\lambda$ satisfy the integrability equations~(\ref{INTEQ}).
Therefore, up to af\/f\/ine symplectic transformation, there exists a
unique $f_{\lambda}$ with $\mathcal{F}_\lambda = \mathcal{F}$. As
$p_\lambda$ are distinct for distinct values of $\lambda$, the
immersions $f_\lambda$ are noncongruent for distinct values of
$\lambda$. This shows that applicable elliptic Lagrangian immersions
come in 1-parameter families.

Interestingly enough, away from the points where $h$
vanishes (umbilics), the Gauss map of an applicable elliptic
Lagrangian immersion $f$ is an isothermic spacelike immersion.
Moreover, the Gauss maps of the 1-parameter family of associates to $f$ are
the $T$-transforms of Bianchi and Calapso of~$\gamma_f$
(cf. \cite{Calapso,Bianchi,HMN}).
\end{remark}

\begin{remark}
The 1-parameter family of associates to an applicable elliptic
Lagrangian immersion are related to its second order af\/f\/ine
symplectic deformations. We will not discuss this issue here. We
just recall that two elliptic Lagrangian immersions $f,\tilde{f} : M
\to \mathbb{R}^4$ are {second order af\/f\/ine symplectic deformations of each
other} if there exists a smooth map ${D} : M \to \mathcal{S}(4,\mathbb{R})$
such that, for each $q\in M$, the Taylor expansions about $q$ of
$D(q)f$ and $\widetilde{f}$ agree through second order terms. For
the general notion of applicability and deformation in homogeneous
spaces we refer the reader to \cite{Gr,JensenJDG}.
\end{remark}

\subsection{The dif\/ferential equations of applicable Lagrangian immersions}

Let $f:M\to \mathbb{R}^4$ be an applicable elliptic Lagrangian immersion and
let $\mathcal{Q}$ be a non-zero holomorphic quadratic dif\/ferential
on $M$ such that $\mathcal{H}\wedge_{\mathbb{R}}\mathcal{Q}=0$. If we
suppose that $\mathcal{Q}$ never vanishes, we can choose a complex
coordinate $(U,z)$ such that $\mathcal{Q} = dz^2$. Let $\mathbf{S}:U\to
\mathcal{S}(4,\mathbb{R})$ be a frame f\/ield along $f$ adapted to $(U,z)$ with
invariant functions $t$, $h$ and $p$. Now, the equation
$\mathcal{H}\wedge_{\mathbb{R}}\mathcal{Q}=0$ implies that $h$ is real
valued.

Next, in analogy with \cite{BHPP, MNhouston} we deduce the
dif\/ferential equations satisf\/ied by the invariant functions of an
elliptic Lagrangian immersion. With respect to the adapted frame
f\/ield $\mathbf{S}$, we have
\[
  \alpha^1_1=-\alpha^2_2=hdx,\qquad \alpha^2_1=\alpha^1_2=-hdy,
\]
where $h$ is a real-valued function. From (\ref{MCA}), we compute
\[
 \left\{\begin{array}{llll}
  \beta^1_1 = (p_1+ h_x + h^2)dx - (p_2 + h_y)dy,\\
    \beta^2_2 = -(p_1- h_x + h^2)dx + (p_2 - h_y)dy,\\
    \beta^2_1 = -p_2dx - (p_1 - h^2)dy.
     \end{array}\right.
\]
Using (\ref{MCB}), we f\/ind
\[
\left\{\begin{array}{llll}
  h_{xy} = -2hp_2,\\
   (p_1)_y = - (p_2)_x - 2(h^2)_y,\\
   (p_1)_x = (p_2)_y + 2(h^2)_x.
    \end{array}\right.
\]
Thus, calculated with respect to a complex parameter $z$ such that
$\mathcal{Q}=dz^2$, the integrability equations of an applicable
elliptic Lagrangian immersion are
\begin{gather}\label{DIFFEQDEF}
 \left\{\begin{array}{llll}
   t_{\bar{z}} = h \bar{t},\\
    h_{xy} = -2hp_2,\\
    \Delta p_2 = - 4(h^2)_{xy}.
     \end{array}\right.
       \end{gather}
If $(t,h,p_2)$ is a solution of this system, the 1-form
\[
  \left[ (p_2)_y + 2 (h^2)_x \right]dx - \left[(p_2)_x + 2(h^2)_y\right]dy
\]
is closed. If $p_1 : M\to \mathbb{R}$ is so that
\[
  dp_1=\left[(p_2)_y + 2 (h^2)_x \right]dx
         -\left[(p_2)_x + 2(h^2)_y\right]dy,
\]
then $t$, $h$ and $p= p_1+ip_2$ give a solution of (\ref{INTEQ}).
Thus, up to congruence, there exists a unique elliptic Lagrangian
immersion $f : M\to \mathbb{R}^4$ with invariant functions $t$, $h$ and $p$.
Dif\/ferent choices of the primitive $p_1$ yield the family of
associate immersions applicable to $f$.

\subsection{Totally umbilic Lagrangian immersions}\label{section4.4}

\begin{definition}
An elliptic Lagrangian immersion $f$ is called {\it totally umbilic}
if its Hopf dif\/fe\-ren\-tial~$\mathcal{H}$ vanishes identically, or,
equivalently, $h$ vanishes identically.
This terminology is based on the observation that $\mathcal{H} \equiv 0$
if and only if the Gauss map $\mathcal{T}_f$ of $f$ is a totally umbilic
spacelike immersion (cf.~\cite{BrDG}).
\end{definition}

\begin{remark}
Note that $f$ is totally umbilic if and only if its Fubini
cubic form $\mathcal{F}$ is holomorphic (cf.~\eqref{2F12}).
In this case, $M$ is either the complex plane or the unit disk.
Therefore, there exists a~global parameter $z$ on $M$ and a~global adapted frame f\/ield $\mathbf{S}:M\to \mathcal{S}(4,\mathbb{R})$ along~$f$.
\end{remark}

\begin{proposition}\label{prop:non-congruentTU}
Let $f,\widetilde{f}:M\to \mathbb{R}^4$ be two noncongruent, totally umbilic
elliptic Lagrangian immersions. Then there exists a biholomorphism
$\Phi:M\to M$ such that $f$ and $\widetilde{f}\circ \Phi$ are
applicable.
\end{proposition}

\begin{proof}
On $M$ there are complex parameters $z, w : M\to W$, $W=\mathbb{C}$ or
$\Delta$, such that $\mathcal{F}=(dz)^3$ and
$\widetilde{\mathcal{F}}=(dw)^3$. Then $\Phi:=w^{-1}\circ z:M\to M$
is a biholomorphic map such that $f$ and $\widetilde{f}\circ \Phi$
have the same Fubini cubic dif\/ferential. This yields the required
result.
\end{proof}

\subsection{Complex curves and Lagrangian immersions}\label{section4.5}

On $\mathbb{R}^4$, consider the complex structure def\/ined by the
identif\/ication of $\mathbb{R}^4$ with $\mathbb{C}^2$ given by
\[
  I: \ {}^t\hskip-1pt  \big(x^1,x^2,x^3,x^4\big)\in \mathbb{R}^4\mapsto {}^t\hskip-1pt  \big(x^1-ix^2,x^3+ix^4\big)\in
   \mathbb{C}^2.
    \]
Then, any complex line generates a Lagrangian plane, which implies that
a complex immersion $f : M\to \mathbb{R}^4 \simeq \mathbb{C}^2$ is automatically Lagrangian.

Let $\mathrm{ISL}(2,\mathbb{C})\cong \mathbb{C}^2\rtimes \mathrm{SL}(2,\mathbb{C})$ be the inhomogeneous
group associated with the {unimodular complex group} $\mathrm{SL}(2,\mathbb{C})$.
The group $\mathrm{ISL}(2,\mathbb{C})$ can be realized as a closed subgroup of
$\mathcal{S}(4,\mathbb{R})$ by
\[
 \left(
  \begin{array}{cc}
   1 & 0 \\
    v & A \\
     \end{array}
      \right)\in \mathrm{ISL}(2,\mathbb{C})
 \mapsto \left(
  \begin{array}{cc}
    1 & 0 \\
     I^{-1}v & [I^{-1}AI] \\
      \end{array}
       \right)\in \mathcal{S}(4,\mathbb{R}),
         \]
where, if we write $A \in \mathrm{SL}(2,\mathbb{C})$ in the form $A = (v^i_j) + i(w^i_j)$,
for $v^i_j$, $w^i_j \in \mathbb{R}$,
\[
 [I^{-1}AI] =
 \left(
  \begin{array}{cccc}
    v^1_1 & w^1_1 &v^1_2 & -w^1_2 \\
    -w^1_1 & v^1_1 & -w^1_2 & -v^1_2 \\
      v^2_1 & w^2_1 & v^2_2 & -w^2_2 \\
       w^2_1 & -v^2_1 & w^2_2 & v^2_2
        \end{array}
         \right).
          \]
The Lie algebra $\mathfrak{isl}(2,\mathbb{C})$ identif\/ies with the
Lie subalgebra of $\mathfrak{s}(4,\mathbb{R})$ consisting of all
$\mathrm{S}(p,\mathbf{x})\in \mathfrak{s}(4,\mathbb{R})$ such that
\[
 \mathbf{x}=\left(
  \begin{array}{cc}
    a & b \\
     c & -{}^t\hskip-1pt  a \\
     \end{array}
      \right),
       \]
where
\begin{gather}\label{sl-valued}
  b,c\in \mathrm{S}(2,\mathbb{R}),\qquad
   \mathrm{tr}\,(b) = \mathrm{tr}\,(c)=0, \qquad  a^1_1-a^2_2=0,\qquad
     a^2_1+a^1_2=0.
      \end{gather}

The precise relation between complex curves and elliptic Lagrangian
immersions is summarized in the following statement.

\begin{proposition}\label{prop:TU-complexcurve}
If $f : M \to \mathbb{C}^2$ is a complex curve without flex
points\footnote{A point $q\in M$ is a f\/lex point if ${f_z}|_q
\wedge {f_{zz}}|_q = 0$, $z$ complex coordinate on $M$.}, then $f$ is
a totally umbilic elliptic Lagrangian immersion. Conversely, let $f
: M\to \mathbb{R}^4$ be a totally umbilic elliptic Lagrangian immersion.
Then there exists an affine symplectic transformation
$\mathcal{D}\in \mathcal{S}(4,\mathbb{R})$ such that $\mathcal{D}\cdot f$ is a
complex curve without flex points.
\end{proposition}

\begin{proof}
Let $f$ be a complex curve without f\/lex points. We already observed
that $f$ is Lagrangian. Moreover, since $f$ has no f\/lex points, it
is also elliptic. To prove that $f$ is totally umbilic, we have just
to observe that the Maurer--Cartan form of an adapted frame along
$f$ with respect to a given complex parameter takes values in
$\mathfrak{isl}(2,\mathbb{C})$. The def\/ining conditions
\eqref{adapted-frame} for an adapted frame, combined with
\eqref{sl-valued}, yield that the invariant $h$ vanishes
identically, and hence $f$ is totally umbilic.

Conversely, let $f$ be a totally umbilic elliptic Lagrangian immersion.
For a f\/ixed complex parameter $z : M\to \mathbb{C}$, let $\mathbf{S}:M\to \mathcal{S}(4,\mathbb{R})$ be
an adapted frame f\/ield along $f$. Since $f$ is totally umbilic, using
\eqref{t&h}, \eqref{upsilon} and \eqref{p}, we f\/ind that
\[
  \tau^1 = t_1dx-t_2dy, \qquad \tau^2=-t_2dx-t_1dy
\]
and
\[
  \gamma=\left(
    \begin{array}{cc}
      dx & dy \\
       dy & -dx \\
        \end{array}
         \right),\qquad
      \alpha=0,\qquad
       \beta=\left(
         \begin{array}{cc}
           p_1dx-p_2dy & -p_2dx-p_1dy \\
             -p_2dx-p_1dy &  -p_1dx+p_2dy \\
              \end{array}\right).
\]
This implies that $\mathbf{S}^{-1}d\mathbf{S}$ is an $\mathfrak{isl}(2,\mathbb{C})$-valued
1-form of bidegree $(1,0)$, which implies that there exists an
element $\mathcal{D} \in \mathcal{S}(4,\mathbb{R})$ such that $\mathcal{D}\mathbf{S}$ is
holomorphic and takes values in $\mathrm{ISL}(2,\mathbb{C})$. Therefore,
$\mathcal{D}f$ is a complex curve without f\/lex points.
\end{proof}

\begin{remark}
Let $f : M \to \mathbb{C}^2$ be a complex curve. Possibly replacing $f$ with
$\mathcal{D}f$, for some $\mathcal{D} \in\mathcal{S}(4,\mathbb{R})$, the homogeneous
part $\mathbf{X} : M\to \mathrm{SL}(2,\mathbb{C})$ of the frame f\/ield
adapted to a complex coordinate $z$ satisf\/ies
\[
 \mathbf{X}^{-1}d\mathbf{X}=\left(
  \begin{array}{cc}
   0 & -ip \\
    1 & 0 \\
     \end{array}
       \right)dz.
\]
In particular, $f$ can be obtained
by integrating the f\/irst column vector of $\mathbf{X}$, i.e.\
\[
  f=\int \mathrm{X}_1 dz.
\]
Note that $\mathbf{X}:M\to \mathrm{SL}(2,\mathbb{C})$ is a
contact complex curve in $\mathrm{SL}(2,\mathbb{C})$. This shows that
totally umbilic Lagrangian surfaces are strictly related to the
geometry of f\/lat fronts in hyperbolic 3-space~\cite{GMM,KUY}. For a f\/ixed choice of the complex parameter,
let $\mathbf{X}_{\lambda}:M\to \mathrm{SL}(2,\mathbb{C})$ be so that
\[
 \mathbf{X}_{\lambda}^{-1}d\mathbf{X}_{\lambda}=\left(
  \begin{array}{cc}
    0 & -i(p-\lambda) \\
      1 & 0 \\
       \end{array}
        \right)dz.
\]
Then, the 1-parameter family of complex curves $f_{\lambda}$
obtained by integrating the f\/irst column vector of $\mathbf{X}_{\lambda}$
are not congruent to each other and are all applicable over $f$.
\end{remark}

\section{Examples}\label{s:examples}

The simplest non-umbilic solutions of (\ref{DIFFEQDEF}) are given by
taking $h$ and $p$ real constants. Possibly rescaling the complex
parameter $z$ by a real constant, we may assume $h=1$. The only
equation to be solved is then
\[
  t_{\bar{z}}= \bar t .
\]
The homogeneous part of the Maurer--Cartan form of the adapted frame
$\mathbf{S}$ is
\[
 \theta=\mathcal{A}dx+\mathcal{B}dy,
\]
where $\mathcal{A}$, $\mathcal{B}\in \mathfrak{sp}(4,\mathbb{R})$ are given by
\[
 \mathcal{A}=\left(
  \begin{array}{cccc}
    1 & 0 & p+1 & 0 \\
     0 & -1 & 0 & -(p+1) \\
      1 & 0 & -1 & 0 \\
      0 & -1 & 0 & 1 \\
       \end{array}
        \right),
\qquad
 \mathcal{B}=\left(
  \begin{array}{cccc}
   0 & -1 & 0 & 1-p \\
    -1 & 0 & 1-p & 0 \\
     0 & 1 & 0 & 1 \\
       1 & 0 & 1 & 0 \\
      \end{array}
       \right).
\]
Since $\mathcal{A}$ and $\mathcal{B}$ commute, the general solution
of $\mathbf{X}^{-1}d\mathbf{X}$ is
\[
 \mathbf{X}(x,y)=\mathbf{Y}\mathrm{Exp}(x\mathcal{A}+y\mathcal{B}),
\]
where $\mathbf{Y}\in \mathrm{Sp}(4,\mathbb{R})$. Without loss of
generality, we may assume $\mathbf{Y}=I_4$. Then the f\/irst two column
vectors of $\mathbf{X}$ are
\[
 \mathrm{X}_1=\left(
          \begin{array}{c}
            \cosh(\sqrt{2-p}y)\left(\cosh(\sqrt{2+p}x)
       +\frac{\sinh(\sqrt{2+p}x)}{\sqrt{2+p}}\right) \\
           -\frac{(\sqrt{p+2}\cosh(\sqrt{p+2}x)
         +p\sinh(\sqrt{p+2}x))\sinh(\sqrt{2-p}y)}{\sqrt{4-p^2}} \\
            \frac{\cosh(\sqrt{2-p}y)\sinh(\sqrt{p+2}x)}{\sqrt{p+2}} \\
             \frac{(\sqrt{p+2}\cosh(\sqrt{p+2}x)
           +2\sinh(\sqrt{p+2}x))\sinh(\sqrt{2-p}y)}{\sqrt{4-p^2}}
             \end{array}
              \right)
\]
and
\[
 \mathrm{X}_2=\left(
          \begin{array}{c}
            \frac{(-\sqrt{p+2}\cosh(\sqrt{p+2}x)
        +p\sinh(\sqrt{p+2}x))\sinh(\sqrt{2-p}y)}{\sqrt{4-p^2}} \\
            \cosh(\sqrt{2-p}y)\left(\cosh(\sqrt{2+p}x)
         -\frac{\sinh(\sqrt{2+p}x)}{\sqrt{2+p}}\right) \\
             \frac{(\sqrt{p+2}\cosh(\sqrt{p+2}x)
            -2\sinh(\sqrt{p+2}x))\sinh(\sqrt{2-p}y)}{\sqrt{4-p^2}} \\
            - \frac{\cosh(\sqrt{2-p}y)\sinh(\sqrt{p+2}x)}{\sqrt{p+2}}
          \end{array}
        \right).
\]
If $t=t_1+it_2$ is a solution of $t_{\bar{z}}=\bar{t}$, that is
\[
  (t_1)_x- (t_2)_y=t_1,\qquad  (t_1)_y + (t_1)_x=-t_2,
\]
then the $\mathbb{R}^4$-valued 1-form
\[
  \mathrm{Y}_1dx+\mathrm{Y}_2dy=(t_1\mathrm{X}_1-t_2\mathrm{X}_2)dx-(t_2\mathrm{X}_1+t_1\mathrm{X}_2)dy
\]
is closed and the corresponding Lagrangian immersion $f$ can be computed
by solving
\[
  df=\mathrm{Y}_1dx+\mathrm{Y}_2dy=(t_1\mathrm{X}_1-t_2\mathrm{X}_2)dx-(t_2\mathrm{X}_1+t_1\mathrm{X}_2)dy.
\]
 This is an elliptic Lagrangian immersion with Fubini's cubic
form $\mathcal{F}=t^2\,dz^3$. For a f\/ixed $t$ and $p\neq 2$, we have
the 1-parameter family of applicable immersions. Explicit solutions
of $t_{\bar{z}}=\bar {t}$ can be obtained with the ansatz
\[
  t_1(x,y)=v_1(x)+w_1(y),\qquad t_2(x,y)=v_2(x)+w_2(y).
\]
We then have
\begin{alignat*}{3}
&  v_1(x)  =c_1e^{2x}-a_1,\qquad & & v_2(x)  =c_2e^{-2x}-a_2, &\\
&   w_1(y)  =m_1e^{2y}+m_2e^{-2y}-a_1, \qquad && w_2(y)  =-m_1e^{2y}+m_2e^{-2y}+a_2,&
\end{alignat*}
where $a_1$, $a_2$, $c_1$, $c_2$, $m_1$ and $m_2$ are real
constants. In this case, the integration involves elementary
functions and can be performed explicitly. The general formulae are
rather involved expressions and can be obtained with any standard
software of scientif\/ic computation such as {\sc Mathematica}~6.

For instance, if $a_1=a_2=0$ and $m_1=m_2=0$, the components of the
immersions $f$ are given by
\begin{gather*}
 f_j(x,y) =\frac{e^{-2x}}{(p-2)\sqrt{(p+2)(4-p^2)}}
   \left(\cosh(\sqrt{p+2}x)A_j(x,y) +\sinh(\sqrt{p+2}x)B_j(x,y)\right),
    \end{gather*}
$j=1,\dots,4$, where the functions $A_j(y)$ and $B_j(y)$ are given by
\begin{gather*}
    A_1(x,y) = - c_1e^{4x}\sqrt{(p+2)(4-p^2)}\cosh(\sqrt{2-p}y)
  - c_2(p^2-4)\sinh(\sqrt{2-p}y),\\
      A_2(x,y)=- c_1e^{4x}(p^2-4)\sinh(\sqrt{2-p}y)
 -c_2\sqrt{(p+2)(4-p^2)}\cosh(\sqrt{2-p}y)),\\
         A_3(x,y)=c_1e^{4x}\sqrt{(p+2)(4-p^2)}\cosh(\sqrt{2-p}y),\\
           A_4(x,y)=c_2\sqrt{(p+2)(4-p^2)}\cosh(\sqrt{2-p}y),
             \end{gather*}
and by
\begin{gather*}
   B_1(x,y)=c_1e^{4x}p\sqrt{4-p^2}\cosh(\sqrt{2-p}y)
    -c_2(p-2)\sqrt{p+2}\sinh(\sqrt{2-p}y),\\
      B_2(x,y)=c_1e^{4x}(p-2)\sqrt{p+2}\sinh(\sqrt{2-p}y)
 -c_2p\sqrt{4-p^2}\cosh(\sqrt{2-p}y),\\
          B_3(x,y)=-2c_1e^{4x}\sqrt{4-p^2}\cosh(\sqrt{2-p}y)
-c_2(p-2)\sqrt{p+2}\sinh(\sqrt{2-p}y),\\
           B_4(x,y)=c_1e^{4x}(p-2)\sqrt{p+2}\sinh(\sqrt{2-p}y) +2c_2\sqrt{4-p^2}\cosh(\sqrt{2-p}y).
              \end{gather*}
These provide explicit examples of elliptic Lagrangian immersions with
the same Gauss map but which are not applicable
and
of applicable non-umbilic elliptic Lagrangian immersions.
Note that for $p>2$ some
of the hyperbolic functions transform into trigonometric functions.

\subsection*{Acknowledgements}

The work was partially supported by MIUR projects: \textit{Metriche
riemanniane e variet\`a dif\-fe\-ren\-zia\-bi\-li} (E.M.);
\textit{Propriet\`a geometriche delle variet\`a reali e complesse}
(L.N.); and by the GNSAGA of INDAM.  The authors would like to thank
the referees for their useful comments and suggestions.

\pdfbookmark[1]{References}{ref}
\LastPageEnding

\end{document}